\documentclass[a4paper,12pt]{article}
\usepackage[english]{babel}
\usepackage{amsmath}
\usepackage{amssymb}
\usepackage{tabls}
\usepackage[pdftex]{graphicx}
\usepackage{calc}
\usepackage{xcolor}
\usepackage[T1]{fontenc}
\usepackage[latin1]{inputenc}
\title{Maximal Graded Orders over Crystalline Graded Rings}
\author{Tim Neijens\\University of Antwerp\\ \texttt{tim.neijens@gmail.com} \and Freddy Van Oystaeyen \\ University of Antwerp \\ \texttt{fred.vanoystaeyen@ua.ac.be}}

\newcommand{\N}{\mathbb{N}}
\newcommand{\Z}{\mathbb{Z}}

\newcommand{\C}{\mathbb{C}}
\newcommand{\blok}{\hfill \Box}
\newcommand{\CGR}{\mathop \diamondsuit \limits_{\sigma ,\alpha}}

\newtheorem{defi}{Definition} 
\newtheorem{opm}[defi]{Remark}
\newtheorem{prop}[defi]{Proposition}

\newtheorem{st}[defi]{Theorem}
\newtheorem{lem}[defi]{Lemma}

\begin{document}
\maketitle

\begin{abstract}
Crystalline graded rings are generalizations of certain classes of rings like generalized twisted group rings, generalized Weyl algebras, and generalized skew crossed products.  When the base ring is a commutative Dedekind domain, two constructions are given for producing maximal graded orders.  On the way, a new concept appears, so-called, spectrally twisted group.  Some general properties of it are studied.  At the end of the paper several examples are considered.
\end{abstract}

\section{Preliminaries}
{\defi \label{def1}\textbf{Pre-Crystalline Graded Ring}\\
Let $A$ be an associative ring with unit $1_A$.  Let $G$ be an arbitrary group.  Consider an injection $u: G \rightarrow A$ with $u_e = 1_A$, where $e$ is the neutral element of $G$ and $u_g \neq 0$,  $\forall g \in G$.  Let $R \subset A$ be an associative ring with $1_{R}=1_A$.  We consider the following properties:
\begin{description} 
	\item[(C1)]\label{def2} $A = \bigoplus_{g \in G} R u_g$.
	\item[(C2)]\label{def3} $\forall g \in G$, $R u_g = u_g R$ and this is a free left $R$-module of rank $1$.
	\item[(C3)]\label{def4} The direct sum $A = \bigoplus_{g \in G} R u_g$ turns $A$ into a $G$-graded ring with $R = A_e$.
\end{description}
We call a ring $A$ fulfilling these properties a \textbf{Pre-Crystalline Graded Ring}.}\\

{\prop \label{def5} With conventions and notation as in Definition \ref{def1}:
\begin{enumerate}
	\item For every $g \in G$, there is a set map $\sigma_g : R \rightarrow R$   defined by: $u_g r = \sigma_g(r)u_g$ for $r \in R$.  The map $\sigma_g$ is in fact a surjective ring morphism.  Moreover, $\sigma_e = \textup{Id}_{R}$.
	\item There is a set map $\alpha : G \times G \rightarrow R$ defined by $u_g u_h = \alpha(g,h)u_{gh}$ for $g,h \in G$.  For any triple $g,h,t \in G$ the following equalities hold:
		\begin{eqnarray}
		\alpha(g,h)\alpha(gh,t)&=&\sigma_g(\alpha(h,t))\alpha(g,ht) \label{def6},\\
		\sigma_g(\sigma_h(r))\alpha(g,h)&=& \alpha(g,h)\sigma_{gh}(r) \label{def7}.
		\end{eqnarray}
	\item $\forall g \in G$ we have the equalities $\alpha(g,e) = \alpha(e,g) = 1$ and $\alpha(g,g^{-1}) = \sigma_g(\alpha(g^{-1},g)).$
\end{enumerate}
}
\begin{flushleft}\textbf{Proof}\end{flushleft} See \cite{NVO6}. $\blok$\\

{\prop Notation as above, the following are equivalent:
\begin{enumerate}
	\item $R$ is $S(G)$-torsionfree.
	\item $A$ is $S(G)$-torsionfree.
	\item $\alpha(g,g^{-1})r=0$ for some $g \in G$ implies $r = 0$.
	\item $\alpha(g,h)r=0$ for some $g,h \in G$ implies $r = 0$.
	\item $R u_g = u_g R$ is also free as a right $R$-module with basis $u_g$ for every $g \in G$.
	\item for every $g \in G$, $\sigma_g$ is bijective hence a ring automorphism of $R$.
\end{enumerate}
}
\begin{flushleft}\textbf{Proof}\end{flushleft} See \cite{NVO6}. $\blok$\\

{\defi Any $G$-graded ring $A$ with properties \textbf{(C1),(C2),(C3)}, and which is $G(S)$-torsionfree is called a \textbf{crystalline graded ring}.  In case $\alpha(g,h) \in Z(R)$, or equivalently $\sigma_{gh}=\sigma_g \sigma_h$, for all $g,h \in G$, then we say that $A$ is \textbf{centrally crystalline}.}\\

In the study of generalized twisted group rings an important role is played by
\[H=\{h \in G | \alpha(h,h^{-1}) \textrm{ is invertible in }R\}.\]
Observe that if $h \in G$ also $\alpha(h^{-1},h)$ is invertible in $R$.\\

{\prop \label{def12} For $H \subset G$ as defined above and $A = R \CGR G$ a CGR, the following properties hold:
\begin{enumerate}
	\item For $h \in H$ and $x \in G$, $\alpha(x,h)$ and $\alpha(h,x)$ are invertible in $R$.
	\item $H$ is a subgroup of $G$.
	\item For $h, h' \in H$ and $x,y \in G$:
	\[R \alpha(hx,yh') = R \sigma_h(\alpha(x,y)).\]
	In particular, if $\alpha(x,y)$ is invertible in $R$ then so is $\alpha(hx, yh')$ for every $h,h' \in H$.
\end{enumerate}
\ \\}
\textbf{Proof}
\begin{enumerate}
	\item For $h \in H$ and $x \in G$ we have the $2$-cocycle relation:
	\[\alpha(x,h)\alpha(xh, h^{-1})=\sigma_x(\alpha(h,h^{-1}))\alpha(x,e).\]
	The right hand side is invertible in $R$, hence $\alpha(x,h)$ is invertible too because the right inverse yields a left inverse ($\alpha(x,h)$ is normalizing in $R$) and obviously these two will then have to coincide.
	On the other hand, the $2$-cocycle relation (\ref{def6}) for $(h^{-1}, h, x)$ yields
	\[\alpha(h^{-1}, h)\alpha(e,x) = \sigma_{h^{-1}}(\alpha(h,x))\alpha(h^{-1}, hx),\]
	where $\alpha(h^{-1},h)$ is invertible in $R$ because $h \in H$, and $\alpha(h^{-1}, hx)$ is invertible in $R$ by the foregoing since $h^{-1} \in H$.  Consequently $\sigma_{h^{-1}}(\alpha(h,x))$ and therefore also $\alpha(h,x)$ is invertible in $R$.
	\item Consider $x,y \in H$ and look at the $2$-cocycle relation (\ref{def6}) for $(x, y, y^{-1}x^{-1})$, i.e. :
	\[\alpha(x,y)\alpha(xy, y^{-1}x^{-1}) = \sigma_x(\alpha(y, y^{-1}x^{-1}))\alpha(x,x^{-1}),\]
	where $\alpha(x,x^{-1})$ and $\sigma_x(\alpha(y, y^{-1}x^{-1}))$ are invertible in $R$ because $x,y \in H$ (using 1.).  It follows that $\alpha(xy, y^{-1}x^{-1})$ is invertible in $R$ too.
	\item From 1. it follows that $R u_{hx} = R u_h R u_x$, for $h \in H$.  Then:
	\[R u_{hx}u_{yh'} = R\alpha(hx,yh')u_{hxyh'}.\]
	But on the other hand
	\begin{eqnarray*}
	R u_{hx}u_{yh'} &=& R u_h u_x R u_y u_{h'}\\
	\ &=& R u_h u_x u_y u_{h'} \\
	\ &=& R u_h \alpha(x,y)u_{xy}u_{h'}\\
	\ &=& R \sigma_h(\alpha(x,y))u_h u_{xy} u_{h'}\\
	\ &=& R \sigma_h(\alpha(x,y))u_h u_{xyh'} = R \sigma_h(\alpha(x,y))u_{hxyh'}.
	\end{eqnarray*} 
	applying the remark starting the proof of 3.  Since $R u_{hxyh'}$ is free (on the left) it follows that $R \alpha(hx,yh')= R\sigma_h(\alpha(x,y))$.  Finally $\alpha(x,y)$ is invertible in $R$ if and only if $\sigma_h(\alpha(x,y))$ is invertible if and only if $R = R \alpha(hx, yh')$, i.e. when $\alpha(hx,yh')$ is invertible. $\blok$
\end{enumerate}

Suitable conditions relating $H$ and $G$ play an important role in the algebraic structure theory of generalized twisted group rings (\cite{NVO3}, \cite{NVO4}).  For example, over a Dedekind domain $D$ the condition $[G:H]<\infty$ yields $G' \subset H$ ($G'$ denoting the commutator subgroup) and $D \alpha(x,y) = D \alpha(y,x)$ then follows from the assumption that these are semiprime ideals.  Further relations between $\alpha$-regular elements in $H$ and $G$ allow to control the center(s) and the Azumaya algebra properties of the corresponding graded rings.\\

\section{Spectrally Twisted Groups}\label{grmax1}

In this section we will introduce an object closely related to a group.  This object, the spectrally twisted group (STG) has more than one operation on it, and those operations are not necessarily associative, but they are connected.

\subsection{Some Definitions, Conventions and Properties}

{\defi \textbf{Spectrally Twisted Group}\\
Consider $\{\Gamma, \theta, G, F, \{M_P\}_{P \in F}\}$.  $\Gamma$ and $F$ are sets.  $G,\cdot$ is a group, acting on the right on $F$.  $\theta$ is a surjective map from $\Gamma$ to $G$.  $\{M_P\}_{P \in F}$ are maps from $\Gamma \times \Gamma$ to $\Gamma$.  They satisfy $\forall P \in F$, $\forall x,y,z \in \Gamma$:
\[M_P[M_P(x,y),z] = M_P[x, M_{P\theta(x)}(y,z)],\]
\[\theta[M_P(x,y)]=\theta(x)\cdot\theta(y).\]
For ease of notation, we will write $\forall P \in F$, $\forall x,y \in \Gamma$:
\[\theta_x := \theta(x),\]
\[(x,y)_P := M_P(x,y).\]
$\exists e_\Gamma \in \Gamma: \forall P \in F, \forall x \in \Gamma: (x,e_\Gamma)_P = x = (e_\Gamma,x)_P$.\\
Given a $P \in F$, $x \in \Gamma$, $\exists {}^Px, x^P \in \Gamma$, the left resp. right inverse of $x$, which satisfy
\[({}^Px,x)_P = e_\Gamma = (x,x^P).\]
We call $\{\Gamma, \theta, G, F, \{M_P\}_{P \in F}\}$ a \textbf{spectrally twisted group}.  If there is no confusion possible, we only write $\Gamma$ to denote the STG.}\\

{\lem \label{grmax2} With notations as above, we have $\forall x \in \Gamma, \forall P \in F$:
\begin{enumerate}
	\item $e_\Gamma$ is unique.
	\item If $y$ is a left or right inverse to $x$ for $P$, then $\theta_y=\theta_x^{-1}$.
	\item ${}^Px = x^{P\theta_x^{-1}}$.
	\item ${}^Px$ and $x^P$ are unique.
\end{enumerate}}
\begin{flushleft}\textbf{Proof}\end{flushleft}
\begin{enumerate}
	\item Suppose $e'_\Gamma$ is also an identity element, then $\forall P \in F$:
	\[e'_\Gamma = (e'_\Gamma, e_\Gamma)_P = e_\Gamma.\]
	\item $(y,x)_P = e_\Gamma \Rightarrow \theta_y \cdot \theta_x = \theta_{e_\Gamma} = e_G \Rightarrow \theta_y = \theta_x^{-1}$.
	\item Fix $x \in \Gamma$ and $P \in F$:
	\begin{eqnarray*}
	{}^Px &=& ({}^Px, e_\Gamma)_P = [{}^Px, (x,x^{P\theta_x^{-1}})_{P\theta_x^{-1}}]_P\\
	\ &=& [({}^Px,x)_P, x^{P\theta_x^{-1}}]_P = x^{P\theta_x^{-1}}.
	\end{eqnarray*}
	\item Suppose $(y,x)_P = e_\Gamma$ then
	\begin{eqnarray*}
	y &=& (y,e_\Gamma)_P\\
	\ &=& [y,({}^Px,x)_P]_P\\
	\ &=& [y,(x,x^{P\theta_x^{-1}})_{P\theta_x^{-1}}]_P\\
	\ &=& [(y,x)_P, x^{P\theta_x^{-1}}]_P\\
	\ &=& x^{P\theta_x^{-1}}\\
	\ &=& {}^Px.
	\end{eqnarray*}
	Similarly for a right inverse.$\blok$
\end{enumerate}

{\opm A group $\kappa, \cdot$ is also an STG if we set $G = \kappa, \theta = \textup{Id}_\kappa, F = {P}, M_P = "\cdot"$.}\\

\subsection{Morphisms}
{\defi Consider two STG's
\[\{\Gamma, \theta, G, F, \{M_P\}_{P \in F}\},\]
\[\{\Gamma', \theta', G', F', \{M'_P\}_{P \in F'}\}.\]
Then $f = \{f_\Gamma, f_G, f_F\}$ is a \textbf{morphism of STG's} if
\begin{enumerate}
	\item $f_\Gamma: \Gamma \rightarrow \Gamma'$ satisfies
	\[f_\Gamma\left(M_P(x,y)\right)=M'_{f_F(P)}\left(f_\Gamma(x), f_\Gamma(y)\right).\]
	\item $f_G:G\rightarrow G'$ satisfies
	\[f_G\left(\theta(x)\cdot \theta(y)\right)=\theta'\left(f_G(x)\right)\cdot \theta'\left(f_G(y)\right).\]
	\item $f_F:F \rightarrow F'$ satisfies
	\[f_F(P\theta_x)=f_F(P)\theta'_{f_G(x)}.\]
\end{enumerate}
For ease of notation, we use $f$ without subscript, even if we mean $f_\Gamma, f_G, f_F$.}\\

{\defi Let $f$ be a morphism, then
\[\textup{Ker f}:= \{x \in \Gamma | f_\Gamma(x)=e_{\Gamma'}\}.\]}\\

In the definition, $\theta$ is a map from $\Gamma$ to its twisting group $G$.  But $\theta$ completely destroys the spectral twist.  We introduce a family of morphisms that retain a bit of that structure:

{\defi Let $\Gamma$ be an STG, and let $T, \cdot$ be a group.  Consider the family $\psi:=\{\psi_P\}_{P \in F}$ of maps from $\Gamma$ to $T$.  Suppose that $\forall P \in F, \forall x,y \in \Gamma$ we have
\[\psi_P\left((x,y)_P\right)= \psi_P(x)\cdot\psi_{P \theta_x}(y),\]
then we call the family $\{\psi_P\}_{P \in F}$ a \textbf{spectrally twisted multiplicative map} (STMM).}\\

{\opm This definition does not conflict with the spectral twist of the $\{M_P\}_{P \in F}$, since for $P \in F, x,y,z \in \Gamma$:
\begin{align*}
\psi_P\left[\left((x,y)_P,z\right)_P\right] &= \psi_P\left[(x,y)_P\right]\psi_{P\theta_x\theta_y}(z)\\
&= \psi_P(x)\psi_{P\theta_x}(y)\psi_{P\theta_x\theta_y}(z)\\
&= \psi_P(x)\psi_{P\theta_x}\left[(y,z)_{P\theta_x}\right]\\
&= \psi_P\left[x, (y,z)_{P\theta_x}\right].
\end{align*}}\\

{\lem \label{grmax3}With notations as above, we have $\forall x \in \Gamma, \forall P \in F$:
\begin{enumerate}
	\item $\psi_P(e_\Gamma)=e_T$.
	\item $\psi_{P\theta_x}(x^P) = \left(\psi_P(x)\right)^{-1}$.
	\item $\psi_{P}({}^Px)=\left(\psi_{P\theta_x^{-1}}(x)\right)^{-1}$.
\end{enumerate}}
\begin{flushleft}\textbf{Proof}\end{flushleft} Straightforward calculation. $\blok$\\

{\defi Let $\Gamma$ be an SGT, $T$ a group. $\psi = \{\psi_P\}_{P\in F}$ an STMM from $\Gamma$ to $T$, then we define
\[\textup{Ker}\psi_P := \{x \in \Gamma | \psi_P(x)=e_T\}.\]}\\

\subsection{Subgroups}
Let $\{\Gamma, \theta, G, F, \{M_P\}_{P \in F}\}$ be an STG.  Suppose $\chi \subset \Gamma$ is closed under $\{M_P\}_{P \in F}$, $e_\Gamma \in \chi$ and $\forall x \in \chi, \forall P \in F: {}^Px, x^P \in \chi$.  It is easy to prove that $\theta(\chi)$ is a subgroup of $G$ in this case.

{\defi Considering the above, we say $\{\chi, \theta_{|\chi}, \theta(\chi), F, \{M_P\}_{P \in F}\}$ is a \textbf{subgroup} of $\{\Gamma, \theta, G, F, \{M_P\}_{P \in F}\}$.  We write $\chi \leq \Gamma$.}\\

\section{Maximal Graded Orders}\label{grmax4}
\subsection{Spectrally Twisted Cocycles}

Let $D$ be a Dedekind domain, $\alpha : G\times G \Rightarrow D \backslash \{0\}$ be a generalized twisted $2$-cocycle, i.e.
\[\alpha(g,h)\alpha(gh,t)= \sigma_g(\alpha(h,t))\alpha(g,ht), \ \ \ \ \forall g,h,t \in G.\]
Define 
\[k : \textup{Spec}D \times G \times G \rightarrow \Z:(P,g,h)\mapsto k_P(g,h),\]
by the decomposition in prime ideals of $D \alpha(g,h)$, i.e.
\[I_{g,h}=D\alpha(g,h)=\prod_{P \in \textup{Spec}D}P^{k_P(g,h)}.\]
We find the following relation for $g,h,t \in G$, $P \in \textup{Spec}D$:
\[D\alpha(g,h)\alpha(gh,t)=D\sigma_g(\alpha(h,t))\alpha(g,ht).\]
\begin{eqnarray*}
\ &\Rightarrow& \prod_P P^{k_P(g,h)} \prod_P P^{k_P(gh,t)}=\prod_P \sigma_g(P)^{k_P(h,t)}\prod_P P^{k_P(g,ht)}\\
\ &\mathop \Rightarrow \limits^{\textup{unicity of decomposition}}& k_P(g,h)+ k_P(gh,t)=k_{Pg}(h,t)+k_P(g,ht).
\end{eqnarray*}

{\defi Let $R$ be a ring, $G$ a group with action on $\textup{Spec}R$, then a map $k: \textup{Spec}R \times G \times G \rightarrow \Z$ is called a \textbf{spectrally twisted $2$-cocycle} if and only if $\forall g,h,t \in G$ and $\forall P \in \textup{Spec}R$:
\[k_P(g,h)+ k_P(gh,t)=k_{Pg}(h,t)+k_P(g,ht).\]}\\

We need to find an equivalence relation between such spectrally twisted $2$-cocycles.  We base this on the equivalence relation for twisted $2$-cocycles as follows.  Let $D$ be a Dedekind domain.  Let $g,h \in G$, $\alpha, \beta:G \times G \rightarrow D$ be twisted generalized $2$-cocycles, $\mu: G \rightarrow D$ a map and $P \in \textup{Spec}D$:
\begin{eqnarray*}
\ &\ & \beta(g,h)\mu(gh)=\alpha(g,h)\mu(g)\sigma_g(\mu(h))\\
\ &\Rightarrow& \prod_P P^{k'_P(g,h)}\prod_P P^{\lambda_P(gh)}= \prod_P P^{k_P(g,h)}\prod_P P^{\lambda_P(g)} \prod_P \sigma_g(P)^{\lambda_P(h)}\\
\ &\Rightarrow& k'_P(g,h)+\lambda_P(gh)=k_P(g,h)+\lambda_P(g)+\lambda_{Pg}(h).
\end{eqnarray*}
This last equation gives us a suitable equivalence relation.\\

{\defi For a group $G$ and a Dedekind domain $D$, let $k$ and $k'$ be spectrally twisted $2$-cocycles.  Then $k$ and $k'$ are said to be \textbf{equivalent $2$-cocycles} if and only if there exists a map $\lambda: \textup{Spec}D \times G \rightarrow \Z$ which satisfies $\forall g,h \in G$:
\[k'_P(g,h)+\lambda_P(gh)=k_P(g,h)+\lambda_P(g)+\lambda_{Pg}(h).\]}\\

Note that since the map $\mu$ giving the equivalence on the level of the twisted $2$-cocycles is not unique, the map $\lambda$ also is not necessarily unique.\\

Now if $|G|=n$, define the set $\widetilde G := n^{-1}\Z \times G$ and $\forall P \in \textup{Spec}D$ operations $M_P$ defined by
\[M_P:\widetilde G \times \widetilde G \rightarrow \widetilde G:[(a,g),(b,h)]\mapsto (a + b + k_P(g,h), gh).\]

{\lem The operation $M_P$ defined as above is not associative but satisfies $\forall a,b,c \in Z, \forall g,h,t \in G$:
\begin{equation}\label{grmax5}M_P\left[M_P\left[(a,g),(b,h)\right],(c,t)\right]=M_P\left[(a,g), M_{Pg}\left[(b,h),(c,t)\right]\right].\end{equation}}\\
\textbf{Proof}
\begin{eqnarray*}
M_P\left[M_P\left[(a,g),(b,h)\right],(c,t)\right] &=& M_P[(a+b+k_P(g,h),gh),(c,t)]\\
\ &=& (a+b+c+k_P(g,h)+k_P(gh,t),ght)\\
\ &=& (a+b+c+k_P(g,ht)+k_{Pg}(h,t),ght)\\
\ &=& M_P[(a,g),(b+c+k_{Pg}(h,t),ht)]\\
\ &=& M_P[(a+b+k_P(g,h),gh),(c,t)].
\end{eqnarray*}
\ $\blok$\\

We can define $\forall P \in \textup{Spec}D$:
\[i: n^{-1}\Z \rightarrow \widetilde G:a\mapsto(a,e),\]
\[\pi: \widetilde G \rightarrow G : (a,g)\mapsto g.\]

{\lem \label{grmax6} The element $(0,e)$ acts as a neutral element for $M_P$, $\forall P \in \textup{Spec}D$.  Let $(a,g)\in \widetilde G$, then the left inverse is $(-a-k_{Pg}(g,g^{-1}),g^{-1})$, and the right inverse is $(-a-k_P(g,g^{-1}),g^{-1})$.}\\
\textbf{Proof} Easy calculation. $\blok$\\

{\prop $\{\widetilde G, \pi, G, \textup{Spec}D, \{M_P\}_{P \in \textup{Spec}D}\}$ is an STG.}\\
\textbf{Proof} Obvious. $\blok$\\

Now let $|G|=n$, then $\alpha^n$ is equivalent to the trivial $2$-cocycle, and by taking powers, it is clear that $nk_P$ as defined above is equivalent to the trivial $2$-cocycle $\forall P \in \textup{Spec}D$.  This means that $\forall P \in \textup{Spec}D$ we find a map $\gamma_P:G \rightarrow \Z$ with
\[nk_P(g,h)= \gamma_P(g)+\gamma_{Pg}(h)-\gamma_P(gh).\]
It is clear that $\gamma_P(e)=0, \ \forall P \in \textup{Spec}D$.  Note that this map $\gamma$ not necessarily is unique!  In the construction below, we fix such a $\gamma$.\\

\begin{flushleft}Consider the direct product $n^{-1}\Z \times G$ and define $\forall P \in \textup{Spec}D$:\end{flushleft}
\[\phi_P : \widetilde G \rightarrow n^{-1}\Z \times G : (a,g)\mapsto (a + n^{-1}\gamma_P(g),g),\]
\[\phi_P\left[\left((a,g),(b,h)\right)_P\right]=\phi_P(a,g)\phi_{Pg}(b,h).\]

We see that $\{\phi_P|P \in \textup{Spec}D\}$ a spectrally twisted multiplicative map.  $\phi_P$ induces a spectrally twisted multiplicative map
\[\psi_P:\widetilde G \rightarrow n^{-1}\Z: (a,g)\mapsto a + n^{-1}\gamma_P(g).\]
It of course has the property $\psi_P \circ i = \textup{Id}$.\\

\begin{flushleft}Fix $P \in \textup{Spec}D$.  Now, given $g \in G$, choose $a_P(g)\in \Z$ such that\end{flushleft}
\[0 \leq \psi_P(a_P(g),g)<1.\]
Now define $\forall P \in \textup{Spec}D$ the $P$-section of $\pi$
\[s_P: G \rightarrow \widetilde G: g \mapsto (a_P(g),g).\]
Since $\psi_P(a,e)=a+0$, $s_P(e)=(0,e)$.
\[m_P: G\times G \rightarrow \Z:i\left(m_P(g,h)\right)=\left(\left(s_P(g),s_{Pg}(h)\right)_P, {}^P\left(s_P(gh)\right)\right)_P.\]
We calculate for $g,h \in G$, $P \in \textup{Spec}D$:
\begin{align*}
&\left(\left(s_P(g),s_{Pg}(h)\right)_P,{}^P\left(s_P(gh)\right)\right)_P\\ &=\left(\left(\left(a_P(g),g\right),\left(a_{Pg}(h),h\right)\right)_P,\left[-a_P\left(gh\right)-k_P\left(gh,(gh)^{-1}\right),(gh)^{-1}\right]\right)_P\\
&=\left(\left[a_P(g) + a_{Pg}(h)+k_P(g,h),gh\right],\left[-a_P(gh)-k_P\left(gh,(gh)^{-1}\right), (gh)^{-1}\right]\right)_P\\
&= \left[a_P(g) + a_{Pg}(h)-a_P(gh)+k_P(gh),e\right].
\end{align*}

\begin{flushleft}So if we define\end{flushleft}
\[m_P: G\times G \rightarrow \Z:i\left(m_P(g,h)\right)=\left(\left(s_P(g),s_{Pg}(h)\right)_P,{}^P\left(s_P(gh)\right)\right)_P,\]
then we find
\[m_P: G\times G \rightarrow \Z:i\left(m_P(g,h)\right)=a_P(g) + a_{Pg}(h)-a_P(gh)+k_P(g,h).\]
This means that $\forall P \in \textup{Spec}D$, $m_P$ and $k_P$ are equivalent $2$-cocycles.

{\defi Let $P \in \textup{Spec}D$, then $m_P$ as defined above is called the \textbf{$P$-maternal $2$-cocycle} of $A = D \CGR G$ corresponding to $\gamma$.  We call $a_P: G \rightarrow \Z$ the \textbf{$P$-maternal power} corresponding to $\gamma$.}\\

\begin{flushleft}For now, we fix $\gamma$.\end{flushleft}

{\lem \label{grmax7} Let $g \in G$, $a \in \Z$, $P\in \textup{Spec}D$, then
\[\psi_{Pg}[(a,g)^{-1}]=-\psi_P[(a,g)].\]}
\begin{flushleft}\textbf{Proof}\end{flushleft}
\begin{eqnarray*}
\psi_{Pg}[(a,g)^{-1}] &=& \psi_{Pg}[-a-k_P(g,g^{-1}), g^{-1}]\\
\ &=& -a-k_P(g,g^{-1}) + n^{-1}\gamma_{Pg}(g^{-1})\\
\ &=& -a-n^{-1}\gamma_P(g)\ \ \ \ \ \ (nk_P(g,g^{-1})=\gamma_P(g)+\gamma_{Pg}(g^{-1}))\\
\ &=& -\psi_P[(a,g)].
\end{eqnarray*}
$\blok$\\

{\prop The spectrally twisted $2$-cocycle $m_P$ as defined above has its values in $\{0,1\}$.}\\
\begin{flushleft}\textbf{Proof}\end{flushleft}
\begin{eqnarray*}
\psi_P\left(i\left(m_P(g,h)\right)\right) &=& \psi_P\left(\left(s_P(g)s_{Pg}(h)\right)s_P(gh)^{-1}\right)\\
\ &=& \psi_P[s_P(g)s_{Pg}(h)]+\psi_{Pgh}[s_P(gh)^{-1}]\\
\ &\mathop = \limits^{(\textup{lemma  }\ref{grmax7})}& \underbrace{\psi_P[s_P(g)]+\psi_{Pg}[s_{Pg}(h)]-\psi_P[s_P(gh)]}_{ - 1 <  \ldots  < 2}.
\end{eqnarray*}
Since by definition $m_P$ has values in $\Z$, the proposition is proven. $\blok$\\

\subsection{Constructing the Maximal Graded Orders}

Let $A = D \CGR G$ be a crystalline graded ring with $D$ a Dedekind domain.  We will now give a few methods in constructing different maximal graded orders.  The first method is by calculating for each possible $\gamma$ the maternal order as defined in the previous section, then doing an additional calculation to find all graded orders, and as such, all maximal graded orders that contain this maternal order.\\
The second method uses a process of conjugation of an order with an element $u_g$ as to obtain a new order.  If we find a maximal graded order, all conjugates of this order are maximal, but as we shall see in the examples, not all contain the ring $A$.

\subsubsection{Graded orders containing unital graded orders.}\label{grmax8}

{\defi Let $T = \bigoplus_{g \in G} I_g u_g$ be a graded order, $I_g = \prod_{P \in \textup{Spec}D} P^{r_P(g)}$.  Then if, $\forall g,h \in G$ and $P \in \textup{Spec}D$:
\[r_P(g) + r_{Pg}(h)-r_P(gh)+k_P(g,h) \in \{0,1\},\]
we call $T$ a \textbf{unital graded order}.}\\

This section will present a calculation that, when given a unital order, we can construct a graded order above it or prove the order in question is maximal.  We will do so by 'making zeroes' i.e. if $r_P(g) + r_{Pg}(h)-r_P(gh)+k_P(g,h) = 1$ for some $g,h \in G$, $P \in \textup{Spec}D$, we will try to modify the $r_P$ such that the new expression equals zero.\\

Consider $\forall P \in \textup{Spec}D$ the respective maternal $2$-cocycles $m_P$ and maternal powers $a_P$ as defined above, corresponding to some fixed $\gamma$.  Now construct $\forall g \in G$ (as above, all $P$ are taken in $\textup{Spec}D$):
\[J_g = \prod_P P^{a_P(g)},\]
and if we take $\{u_g | g \in G\}$ the basis for $A$ over $D$ then define
\[M=\bigoplus_{g\in G}J_g u_g.\]

{\defi $M$ defined as above is called the \textbf{maternal order} corresponding to $\gamma$.  It is by construction unital.}\\

This is a graded order.  Suppose $T = \oplus_{g \in G} I_g u_g$ is a graded order.  Then for each $g\in G$, $I_g$ must be finitely generated, $I_g$ is a fractional ideal.  This is true for $M$ since we only consider finitely many primes $P$ (those that appear in the decomposition of the $\alpha(g,h)$) in the construction of $m_P$.  (Not entirely true, but for the non-appearing primes we just take everything equal to $0$.)  This means that only finitely many $a_P(g)$ differ from $0$.  We need the following relation $\forall g,h \in G$
\begin{eqnarray*}
I_g u_g I_h u_h \subset I_{gh}u_{gh} &\Leftrightarrow& I_g \sigma_g(I_h)u_g u_h \subset I_{gh} u_{gh}\\
\ &\Leftrightarrow& I_g \sigma_g(I_h) \alpha(g,h)u_{gh} \subset I_{gh}u_{gh}\\
\ &\Leftrightarrow& I_{gh}^{-1}I_g\sigma_g{I_h}\alpha(g,h)D \subset D.
\end{eqnarray*}
Now construct $\forall P \in \textup{Spec}D$, $g,h \in G$:
\[t_P(g,h)=-r_P(gh)+r_P(g)+r_{Pg}(h)+k_P(g,h),\]
where we set
\[I_g = \prod_P P^{r_P(g)}.\]
For $T$ to be an order, $t_P(g,h) \geq 0$, $\forall g,h \in G, \forall P \in \textup{Spec}D$.  This is true for the maternal order $A$ by definition.\\

Suppose $T = \oplus_{g \in G}I_g u_g$ is a graded order with $M \subset T$, $M \neq T$.  So we can assume that $\forall P \in \textup{Spec}D$ and all $g \in G$: $r_P(g)\leq a_P(g)$, and at least one $P_0$ prime and $g_0 \in G$ this inequality is strict.  For the remainder we drop the subscript $\cdot_0$ if no confusion can arise.  Since
\begin{eqnarray*}
t_P(g,g^{-1}) &=& r_P(g)+r_{Pg}(g^{-1})+k_P(g,g^{-1})\\
\ &=&\underbrace{r_P(g)-a_P(g)}_{<0} + \underbrace{r_{Pg}(g^{-1})-a_{Pg}(g^{-1})}_{\leq 0}+m_P(g,g^{-1}).
\end{eqnarray*}
This implies, since $t_P(g,g^{-1})\geq 0$:
\begin{equation}\label{grmax9}
\left\{ {\begin{array}{*{20}c}
   {m_P (g,g^{ - 1} ) = 1} \hfill  \\
   {r_{Pg} (g^{ - 1} ) = a_{Pg} (g^{ - 1} )} \hfill  \\
   {a_P (g) = r_P (g) + 1} \hfill  \\
   {t_P (g,g^{ - 1} ) = 0} \hfill  \\

 \end{array} } \right.
\end{equation}

So, to create a graded order $T \supset M$ where $T \neq M$, we must be able to find $P$ prime, $g \in G$ with $m_P(g,g^{-1})=1$, otherwise we cannot satisfy (\ref{grmax9}).  So, suppose we have $m_P(g,g^{-1})=1$ for some $g$, some $P$, and we have $T \supset M$.  Then
\[t_P(g,g^{-1})=\underbrace{r_P(g)-a_P(g)}_{\leq 0}+\underbrace{r_{Pg}(g^{-1})-a_{Pg}(g^{-1})}_{\leq 0}+ \underbrace{m_P(g,g^{-1})}_{=1} \geq 0.\]

We can now take either $r_P(g)-a_P(g)<0$ or $r_{Pg}(g^{-1})-a_{Pg}(g^{-1})<0$.  Without loss of generality (since $m_P(g,g^{-1})=m_{Pg}(g^{-1},g)$), we can suppose:
\begin{eqnarray*}
r_P(g)&=&a_P(g)-1\\
r_{Pg}(g^{-1})&=&a_{Pg}(g^{-1}).
\end{eqnarray*}
But, we need to remark that this assumption might cause problems since $T$ might not be a graded order any more...  We need to investigate whether or not our changes make any of the $t_Q(A,B)$ negative, in which case we are not allowed to do the proposed change.  Let $Q \in \textup{Spec}D$ and $A, B \in G$:\\

\begin{eqnarray}
t_Q(A,B)&=&[-r_Q(AB)+a_Q(AB)]\nonumber\\
&+&[r_Q(A)-a_Q(A)]\nonumber\\ 
&+&[r_{QA}(B)-a_{QA}(B)]+m_Q(A,B).\label{grmax10}
\end{eqnarray}

To construct a maximal graded order $T \supset M$, we need to construct two main sets, $U$ and $V$.  $U$ will contain all the couples $(P,g)$ for which we can set $r_P(g)=a_P(g)-1$.  $V$ will contain the couples $(P,g)$ for which we cannot set $r_P(g)=a_P(g)-1$, or which we always have to set $r_P(g)=a_P(g)$.\\
Add the $(P,g)$ that cannot be modified, i.e. for which $m_P(g,g^{-1}) = 0$ to $V$.  Now pick a couple $(P,g)$ with $m_P(g,g^{-1})=1$.  We will investigate if we can put $(P,g)$ in $U$.  We do so by checking if equation \ref{grmax10} does not give conflicts (i.e. $t_Q(A,B)<0$).\\

To do this, we will focus on a couple $(P,g)$ and whether or not it can be put in $U$.  During the investigation, it is possible that we are faced with new investigations for couples $(Q,h)$.  We construct a set $W_{(P,g)}$ for such couples.  Note that we need to check all elements in $W_{(P,g)}$ in a similar investigation before we can add $(P,g)$ to $U$.\\

\textbf{To start the investigation}, we create sets $U_{(P,g)}$ and $V_{(P,g)}$ that will be copies of $U$ and $V$, and $W_{(P,g)}$ for couples that need further investigation.  These sets will change throughout the calculation.  We then restart the process for all other possible couples not in $V$.  For now, pick $(P,g)$ not in $V$ and put it in $U_{(P,g)}$.  Put $(Pg,g^{-1})$ in $V_{(P,g)}$.\\

We have to (constantly) check if $\forall Q \in \textup{Spec}D$, $\forall A,B \in G$ the equation in (\ref{grmax10}) is positive if we would carry out $U_{(P,g)}$ and $V_{(P,g)}$ assuming $W_{(P,g)}\subset U_{(P,g)}$.  It is obvious that the terms $m_Q(A,B)$ and $[-r_Q(AB)+a_Q(AB)]$ (whether or not they change) cannot make this expression negative.  Changing $r_P(g)$ might cause a problem when
\begin{eqnarray*}
r_{QA}(B) = r_P(g) &\Rightarrow& A=y \in G, B = g, Q = Ph^{-1}, \label{grmax11}\\ 
r_Q(A) = r_P(g) &\Rightarrow& A = g, B = h \in G, Q = P. \label{grmax12}
\end{eqnarray*}

\paragraph{\underline{$A=y \in G, B = g, Q = Ph^{-1}$}}
\ \\\\
We rewrite (\ref{grmax10}):
\begin{eqnarray*}
t_{Ph^{-1}}(h,g)&=& [-r_{Ph^{-1}}(hg)+a_{Ph^{-1}}(hg)]\\
 &+& [r_{Ph^{-1}}(h)-a_{Ph^{-1}}(h)]\\
 &+& [r_P(g)-a_P(g)] +m_{Ph^{-1}}(h,g).
\end{eqnarray*}

There are two cases:
\begin{enumerate}
	\item $m_{Ph^{-1}}(h,g)=1$.
	\item $m_{Ph^{-1}}(h,g)=0$.
\end{enumerate}

\subparagraph{\underline{$m_{Ph^{-1}}(h,g)=1$.}}
\ \\\\
If $r_{Ph^{-1}}(h)-a_{Ph^{-1}}(h)=-1$ ($(Ph^{-1},h)$ is in $U_{(P,g)}$) then we need to add $(Ph^{-1},hg)$ to $U_{(P,g)}$.  This needs a further check, so put $(Ph^{-1},hg)$ in $W_{(P,g)}$.

\subparagraph{\underline{$m_{Ph^{-1}}(h,g) =0$.}}
\ \\\\
We need to add $(Ph^{-1},h)$ to $V_{(P,g)}$.  If $(Ph^{-1},h)$ is in $U_{(P,g)}$, put $(P,g)$ in $V$.  We can restart the investigation with a new $(P,g)$.

We now want to put $(Ph^{-1},hg)$ in $U_{(P,g)}$, and this can only be done if $m_{Ph^{-1}}(hg, (hg)^{-1})=1$.  It will turn out this is always true.  Our assumption leads to:
\[\underbrace{m_{Ph^{-1}}(h,g)}_{=0} + m_{Ph^{-1}}(hg,g^{-1}) = \underbrace{m_P(g,g^{-1})}_{=1}+\underbrace{m_{Ph^{-1}}(h,e)}_{=0}.\]
And so we find
\begin{equation}m_{Ph^{-1}}(hg,g^{-1})=1. \label{grmax13}\end{equation}
We now consider $2$ possibilities:
\begin{eqnarray*}
m_{Ph^{-1}}(h,h^{-1}) &=& m_P(h^{-1},h) = 1 \label{grmax14},\\
m_{Ph^{-1}}(h,h^{-1}) &=& m_P(h^{-1},h) = 0 \label{grmax15}.
\end{eqnarray*}

\subparagraph{\underline{$m_{Ph^{-1}}(h,h^{-1}) = m_P(h^{-1},h) = 1$.}}
\[
\left\{ {\begin{array}{*{20}c}
   {m_P (g,g^{ - 1} ) = m_{Pg} (g^{ - 1} ,g) = 1} \hfill  \\
   {m_{ Ph^{ - 1}} (h,g) = 0} \hfill  \\
   {m_{ Ph^{ - 1}} (h,h^{ - 1} ) = m_{ Ph^{ - 1}} (h^{ - 1} ,h) = 1} \hfill  \\

 \end{array} } \right.
\]
From this we find
\[\underbrace{m_{Ph^{-1}}(hg,g^{-1})}_{(\ref{grmax13}) \Rightarrow =1} + \underbrace{m_{Ph^{-1}}(h,h^{-1})}_{=1}=m_{Pg}(g^{-1},h^{-1})+m_{h^{-1}}(hg,(hg)^{-1}).\]
This means
\[m_{Ph^{-1}}(hg, (hg)^{-1})=1.\]
We need to repeat the investigation with $(Ph^{-1},hg)$ before we can put it in $U_{(P,g)}$.  Put $(Ph^{-1},hg)$ in $W_{(P,g)}$.\\

\subparagraph{\underline{$m_{Ph^{-1}}(h,h^{-1}) = m_P(h^{-1},h) = 0$.}}
\[
\left\{ {\begin{array}{*{20}c}
   {m_P (g,g^{ - 1} ) = m_{Pg} (g^{ - 1} ,g) = 1} \hfill  \\
   {m_{h^{ - 1} P} (h,g) = 0} \hfill  \\
   {m_{h^{ - 1} P} (h,h^{ - 1} ) = m_{h^{ - 1} P} (h^{ - 1} ,h) = 0} \hfill  \\

 \end{array} } \right.
\]
From this we find
\[m_{Pg}((hg)^{-1},h)+m_{Pg}(g^{-1},h^{-1})=\underbrace{m_{Ph^{-1}}(h,h^{-1})}_{=0},\]
and so
\[m_{Pg}(g^{-1},h^{-1})=0.\]
This in time leads to
\[\underbrace{m_{Ph^{-1}}(hg,g^{-1})}_{=1}+\underbrace{m_{Ph^{-1}}(h,h^{-1})}_{=0}=\underbrace{m_{Pg}(g^{-1},h^{-1})}_{=0} + m_{Ph^{-1}}(hg,(hg)^{-1}),\]
to find
\[m_{Ph^{-1}}(hg, (hg)^{-1})=1.\]
We need to repeat the investigation with $(Ph^{-1},hg)$ before we can put it in $U_{(P,g)}$.  Put $(Ph^{-1},hg)$ in $W_{(P,g)}$.\\

\paragraph{\underline{$A = g, B = h \in G, Q = P$.}}
\ \\\\
We rewrite (\ref{grmax10})
\begin{eqnarray*}
t_P(g,h)&=&[-r_P(gh)+a_P(gh)]\\
&+& [r_P(g)-a_P(g)]\\
&+& [r_{Pg}(h)-a_{Pg}(h)]+m_P(g,h).
\end{eqnarray*}

There are two cases:
\begin{enumerate}
	\item $m_P(g,h)=1$.
	\item $m_P(g,h)=0$.
\end{enumerate}

\subparagraph{\underline{$m_P(g,h)=1$.}}
\ \\\\
If $r_{Pg}(h)-a_{Pg}(h)=-1$ ($(Pg,h)$ is in $U_{(P,g)}$) then we need to add $(P,gh)$ to $U_{(P,g)}$.  This needs a further check, so put $(P,gh)$ in $W_{(P,g)}$.

\subparagraph{\underline{$m_P(g,h)=0$.}}
\ \\\\
We need to add $(Pg,h)$ to $V_{(P,g)}$.  If $(Pg,h)$ is in $U_{(P,g)}$, put $(P,g)$ in $V$.  We can restart the investigation with a new $(P,g)$.
We now want to put $(P,gh)$ in $U_{(P,g)}$, and this can only be done if $m_{P}(gh, (gh)^{-1})=1$.  It will turn out this is always true.  Our assumption leads to:
\[\underbrace{m_P(g,h)}_{=0}+m_P(gh,(gh)^{-1})=m_{Pg}(h,(gh)^{-1})+\underbrace{m_P(g,g^{-1})}_{=1}.\]
We need to repeat the investigation with $(Ph^{-1},hg)$ before we can put it in $U_{(P,g)}$.  Put $(Ph^{-1},hg)$ in $W_{(P,g)}$.\\

We can repeat this proces until we have found for all $(P,g)$ not in $V$ sets $U_{(P,g)}$, $V_{(P,g)}$ and $W_{(P,g)}$.  The set
\[\kappa= \left\{(P,g)|P \in \textup{Spec}D, g \in G\right\}\]
is finite.  Given a $(P,g)$, we can look at all $V_{(Q,h)}$ that appear if we write out this tree:

\[
\begin{array}{*{20}c}
   {} & {} & {(P,g)} & {} & {} & {} & {} & {}  \\
   {} &  \swarrow  &  \downarrow  &  \searrow  & {} & {} & {} & {}  \\
   {U_{(P,g)} } & {} & {V_{(P,g)} } & {} & {W_{(P,g)} } & {} & {} & {}  \\
   {} & {} & {} &  \swarrow  &  \downarrow  &  \searrow  & {} & {}  \\
   {} & {} & {(Q_1 ,h_1 )} & {} &  \ldots  & {} & {(Q_s ,h_s )} & {}  \\
   {} &  \swarrow  &  \downarrow  &  \searrow  & {} &  \swarrow  &  \downarrow  &  \searrow   \\
   {} & {} & {} & {} & {} & {} & {} & {}  \\

 \end{array} 
\]

Then consider $V'$ the union of all the $V_{(Q,h)}$ that appear in the full tree, and $U'$ the union of all $U_{(Q,h)}$ and $W_{(Q,h)}$ that appear in the full tree.  If $U' \cap V' \neq \emptyset$, there is a conflict somewhere and we put $(P,g)$ in $V$, and investigate other $(P,g)$.  If $U' \cap V' = \emptyset$ then there are no conflicts whatsoever, and then we add all appearing $W_{(Q,h)}$ in the full tree to $U_{(P,g)}$, and all appearing $V_{(Q,h)}$ in the full tree to $V_{(P,h)}$.  At last, we can add $U_{(P,g)}$ to $U$ and $V_{(P,g)}$ to $V$.  We can repeat the proces with a new $(P,g) \in \kappa \backslash \left(U \cup V\right)$, until $U \cup V = \kappa$.  If we look at the set $\{U_{(P,g)}|(P,g)\in U\}$, it contains a partition of $U$ by construction.  Every element of this partition $U_{(P,g)}$, represented by $(P,g)$ determines one step in the reduction of $M$, corresponding to a chain of graded orders.  If we do all changes suggested by elements of $U$, then we have found the top of our chain, a maximal graded order.

\begin{flushleft}\textbf{An algorithm}\end{flushleft}
From these calculations we derive an algorithm for determining the unital maximal graded orders, given a crystalline graded ring $A = D \CGR G$, where $D$ is a Dedekind domain.
\begin{enumerate}
	\item Determine the relevant primes, namely those appearing in the prime decomposition of $\alpha(g,h), g,h \in G$.
	\item Determine for all the relevant primes $P$, the corresponding spectrally twisted $2$-cocycle $k_P$ from
	\[D\left(\alpha(g,h)\right) = \prod_{P \in \textup{Spec}D}P^{k_P(g,h)}, \qquad g,h \in G.\]
	\item Calculate $\mu$ from (if $n = |G|$):
	\[\alpha^n(g,h)=\mu(gh)^{-1}\mu(g)\sigma_g(\mu(h))\]
	\item Calculate for the relevant primes $P$ from $\mu$ the $\gamma_P$ from
	\[D(\mu(g))= \prod_{P \in \textup{Spec}D}P^{\gamma_P(g)}, \qquad g \in G\]
	\item Calculate for the relevant primes $P$ from $\gamma_P$ the $a_P \in \Z$ from
	\[0 \leq a_P(g) + \frac{1}{n}\gamma_P(g)<1, \qquad g \in G\]
	\item Calculate for the relevant primes $P$ the spectrally twisted $2$-cocycles $m_P$ from
	\[m_P(g,h)= a_P(g)+a_{Pg}(h)-a_P(gh)+k_P(g,h).\]
	\item \label{grmax16} Determine for which relevant primes $P$ and for which values $g \in G$ we have $m_P(g,g^{-1})=1$.
	\item \label{grmax17} Construct $U$, the set of elements of
	\[\kappa= \left\{(P,g)|P \in \textup{Spec}D, g \in G\right\}\]
	which we will modify ($r_P(g)=a_P(g)-1$) and $V$ the set of elements which we will not modify ($r_P(g)=a_P(g)$).
	\item Set $U = \emptyset$, $V = \left\{(P,g)|m_P(g,g^{-1})=0\right\}$.
	\item \label{grmax18} Pick $(P,g)\notin V$.
	\item Construct $U_{(P,g)} = \left\{(P,g)\right\}$, $V_{(P,g)}=V \cup \left\{(Pg,g^{-1})\right\}$, $W_{(P,g)}=\emptyset$.
	\item \label{grmax19} Find all $h \in G$ for which
	\[m_{Ph^{-1}}(h,g)=1.\]
	If $(Ph^{-1},h)$ is in $U_{(P,g)}$ then add $(Ph^{-1},hg)$ to $W_{(P,g)}$.  If $(Ph^{-1},hg) \in V_{(P,g)}$, put $(P,g)$ in $V$.  Restart from step \ref{grmax18}.
	\item \label{grmax20} Find all $h \in G$ for which
	\[m_{Ph^{-1}}(h,g)=0.\]
	Add $(Ph^{-1},h)$ to $V_{(P,g)}$.  If $(Ph^{-1},h)$ is in $U_{(P,g)}$, put $(P,g)$ in $V$.  Restart from step \ref{grmax18}.  Else, add $(Ph^{-1},hg)$ to $W_{(P,g)}$.
	\item Find all $h \in G$ for which
	\[m_P(g,h)=1.\]
	If $(Pg,h)$ is in $U_{(P,g)}$ then add $(P,gh)$ to $W_{(P,g)}$.  If $(P,gh) \in V_{(P,g)}$, put $(P,g)$ in $V$.  Restart from step \ref{grmax18}.
	\item Find all $h \in G$ for which
	\[m_P(g,h)=0.\]
	Add $(Pg,h)$ to $V_{(P,g)}$.  If $(Pg,h)$ is in $U_{(P,g)}$, put $(P,g)$ in $V$.  Restart from step \ref{grmax18}.  Else, add $(P,gh)$ to $W_{(P,g)}$.
	\item Repeat step \ref{grmax18} for all $(Q,h) \in W_{(P,g)}$.  Since $\kappa$ is finite, this loop ends in finite time.
	\item We now have calculated $U_{(Q,h)}, V_{(Q,h)}, W_{(Q,h)}$ for all $(Q,h)$ appearing in the full tree as described above.  Now check for these $(Q,h)$:
	\[\left(\bigcup U_{(Q,h)} \cup \bigcup W_{(Q,h)}\right)\cap \bigcup V_{(Q,h)}.\]
	If this intersection is empty, add all appearing $W_{(Q,h)}$ to $U_{(P,h)}$ and $U_{(P,g)}$ to $U$.  Then add all appearing $V_{(Q,h)}$ to $V$.  Restart from step \ref{grmax18}.  If this intersection is nonempty, add $(P,g)$ to $V$ and restart from step \ref{grmax18}.
	\item For all $(P,g) \in U$ set $r_P(g)=a_P(g)-1$.  For all $(P,g) \in V$ set $r_P(g)=a_P(g)$.
	\item The graded order
	\[T=\bigoplus_{g \in G}\left(\prod_{P \in \textup{Spec}D}P^{r_P(g)}\right)u_g\]
	is a maximal graded order containing
	\[M=\bigoplus_{g \in G}\left(\prod_{P \in \textup{Spec}D}P^{a_P(g)}\right)u_g.\]
\end{enumerate}

{\opm It is clear that choosing $(P,g)$ in step \ref{grmax18} determines the order $T$ in some way.  To find all maximal graded orders over $M$ we need to use this algorithm for all possible choices $(P,g)$ we make in various points of the algorithm.}

{\st \label{grmax24}Let $A = D \CGR G$ be a crystalline graded ring, where $D$ is a commutative Dedekind domain, $G$ a finite group. Then
\begin{enumerate}
	\item A unital order containing $A$ is maternal.
	\item All maximal graded orders containing $A$ are unital.
\end{enumerate}}
\begin{flushleft}\textbf{Proof}\end{flushleft}
\begin{enumerate}
	\item Suppose that $T = \bigoplus_{g \in G} I_g u_g$ is a graded order, $I_g = \prod_{P \in \textup{Spec}D} P^{r_P(g)}$, $A \subset T$  We need to construct a function $\gamma: \textup{Spec}D\times G \rightarrow \N: (P,x) \mapsto \gamma_P(x)$ such that $\forall P \in \textup{Spec}D, \forall x,y \in G$ ($k_P(x,y)$ defined as usual):
	\[n\cdot k_P(x,y)= \gamma_P(x)+\gamma_{Px}(y)-\gamma_P(xy),\]
	and
	\[0 \leq \frac{1}{n}\gamma_P(x)+r_P(x)<1.\]
	Since $T$ is unital, $\forall P \in \textup{Spec}D, \forall x,y \in G$:
	\[t_P(x,y)= k_P(x,y)+r_P(x)+r_{Px}(y)-r_P(xy)\in \{0,1\}.\]
	Define $\forall P \in \textup{Spec}D, \forall x \in G$, $\gamma_P(x)$ by:
	\[\gamma_P(x)=\left(\sum_{z \in G}t_P(x,z)\right)-nr_P(x)\geq 0.\]
	Then $\forall P \in \textup{Spec}D, \forall x, y \in G$:
	\begin{eqnarray*}
	&\ &\gamma_P(x)+\gamma_{Px}(y)-\gamma_P(xy)\\
	&=& \left(\sum_{z \in G}t_P(x,z)+\sum_{z \in G}t_{Px}(y,z)-\sum_{z \in G}t_P(xy,z)\right)\\
	&\ & \ \ \ \ \ \ \ \ \  - n\left(r_P(x)-r_{Px}(y)+r_P(xy)\right)\\
	&=& \sum_{z \in G}\left(t_P(x,yz)+t_{Px}(y,z)-t_P(xy,z)\right)-nt_P(x,y)+nk_P(x,y)\\
	&=& \left(\sum_{z \in G}t_P(x,y)\right)-nt_P(x,y)+nk_P(x,y) = nk_P(x,y).
	\end{eqnarray*}
	Furthermore $\forall P \in \textup{Spec}D, \forall x \in G$:
	\[0 \leq \frac{1}{n}\gamma_P(x)+r_P(x)= \frac{1}{n}\sum_{z \in G}t_P(x,z)<1,\]
	since $\forall P \in \textup{Spec}D, \forall x,y \in G$, $t_P(x,y) \in \{0,1\}$ and $t_P(x,e)=0$.  In other words, $T$ is a maternal order corresponding to $\gamma$.
	\item Suppose that $T \supset A$ is a maximal graded order.  Using the same notation as in 1. we construct $\forall P \in \textup{Spec}D, \forall x \in G$, $\gamma_P(x)$ by:
	\[\gamma_P(x)=\left(\sum_{z \in G}t_P(x,z)\right)-nr_P(x)\geq 0.\]
	Construct the maternal order $M$ corresponding to $\gamma$.  $M = \bigoplus_{g \in G} J_g u_g$, $J_g = \prod_{P \in \textup{Spec}D} P^{a_P(g)}$.  Fix $P \in \textup{Spec}D$ and $x \in G$.  Suppose
	\[0\leq \frac{1}{n}\gamma_P(x)+r_P(x)<1,\]
	then $r_P(x)=a_P(x)$, since an integer with this property is unique.  Suppose
	\[\frac{1}{n}\gamma_P(x)+r_P(x)\geq 1 > \frac{1}{n}\gamma_P(x)+r_P(x),\]
	then $r_P(x)> a_P(x)$.  Combining these two, $\forall P \in \textup{Spec}D, \forall x\in G$: $0\geq r_P(x)\geq a_P(x)$, implying $M \supset T$.  Since $T$ is maximal, $M = T$. $\blok$
	
\end{enumerate}

\subsection{Conjugation Method}

In this section we will construct new orders by conjugating with $u_g$ for some $g \in G$.  We will therefore define a map $\Psi_g$ for $g \in G$.  Consider $A = D \CGR G$ be a crystalline graded ring with $D$ a Dedekind domain, $K$ the field of quotients of $D$.  Let $g \in G$, then we have the following map on $K \CGR G$:
\[\Psi_g: K \CGR G \rightarrow K \CGR G : \sum_{x \in G} a_x u_x \mapsto u_g \sum_{x \in G} a_x u_x u_g^{-1}.\]

{\prop \label{grmax22}Let $g \in G$.  Set (as in Proposition \ref{def12}):
\[H = \{g \in G|\alpha(g,g^{-1})\in U(D)\}.\]
Then, with notation as above:
\begin{enumerate}
	\item $\Psi_g$ is an automorphism on $K \CGR G$ and is $D^g$-linear.  ($\{D^g = \{d \in D| \sigma_g(d)=d\}\}$.)
	\item $\Psi_h \in \textup{Aut}(A), \forall h \in H$.
	\item If $H \triangleleft G$, then $\Psi_g$ is an automorphism on $S = D \CGR H$.
\end{enumerate}
}
\begin{flushleft}\textbf{Proof}\end{flushleft}
\begin{enumerate}
	\item Very easy to verify.
	\item Since $u_h \in U(A)$, this is clear.
	\item An easy calculation yields ($g \in G, h \in H$):
	\[u_g u_h u_g^{-1}= \alpha(g,h)\alpha^{-1}(ghg^{-1},g)u_{ghg^{-1}}.\]
	Since $H \triangleleft G, ghg^{-1}\in H$.  From Proposition $\ref{def12}$:  $\alpha^{-1}(ghg^{-1},g) \in D$ and as such $u_g u_h u_g^{-1} \in S$.  $\blok$
\end{enumerate}
\ \\

Let $A = D \CGR G$ be a crystalline graded ring with $D$ a Dedekind domain.  Consider a graded order $T$ where $T = \bigoplus_{x \in G}I_x u_x$, $I_x = \prod_{P}P^{r_P(x)}$.  We now conjugate with an element $u_g$ for some $g \in G$.\\

\begin{flushleft}In the following proposition, we will use these calculations for $x,y,g \in G$, $P \in \textup{Spec}D$:\end{flushleft}
\begin{description}
\item[formula 1]
	  \begin{align*}
	  \alpha(xg,x^{-1})\alpha(xgx^{-1},x) &= \sigma_{xg}\left(\alpha(x^{-1},x)\right)\\
	  &= \sigma_{xgx^{-1}}\left(\sigma_x\left(\alpha(x^{-1},x)\right)\right)\\
	  &= \sigma_{xgx^{-1}}\left(\alpha(x,x^{-1})\right)\\
	  \Rightarrow& \alpha(xg,x^{-1})\sigma_{xgx^{-1}}\left(\alpha^{-1}(x,x^{-1})\right)=\alpha^{-1}(xgx^{-1},x).
	  \end{align*}
\item[formula 2]  
	  \begin{align*}
    &t_P(g,g^{-1}xg) = r_P(g)+r_{Pg}(g^{-1}xg)-r_P(xg)+k_P(g,g^{-1}xg)\\
    &\Rightarrow r_{Pg}(g^{-1}xg)+k_P(g,g^{-1}xg)=t_P(g,g^{-1}xg)-r_P(g)+r_P(xg).
    \end{align*}
\item[formula 3]   
    \begin{align*}
    &k_P(x,y)+k_P(xy,g)=k_{Px}(y,g)+k_P(x,yg)\\
    &\Rightarrow k_P(x,y)+k_P(xy,g)-k_{Px}(y,g)=k_P(x,yg).
    \end{align*}
\item[formula 4]  
    \begin{align*}
    &t_P(x,g)=r_P(x)+r_{Px}(g)-r_P(xg)+k_P(x,g)\\
    &\Rightarrow -t_P(x,g)+r_P(x)= -r_{Px}(g)+r_P(xg)-k_P(x,g).
    \end{align*}
\item[formula 5]  
    \begin{equation*}t_P(x,yg)=r_P(x)+r_{Px}(yg)-r_P(xyg)+k_P(x,yg).\end{equation*}
\item[formula 6]   
    \begin{equation*}t_P(x,g)+t_P(xg,g^{-1}yg)=t_{Px}(g,g^{-1}yg)+t_P(x,yg).\end{equation*}
\item[formula 7]   
    \begin{equation*}t_P(g,g^{-1}xg)+t_P(xg,g^{-1}yg)=t_{Pg}(g^{-1}xg,g^{-1}yg)+t_P(g,g^{-1}xyg).\end{equation*}
\end{description}
\ \\
{\prop Define $\Psi_g(T)$ to be $u_g T u_g^{-1}$.  Then $\Psi_g(T)=\bigoplus_{x \in G} \tilde I_x u_x$ is an order.  If we set for $x,y \in G$ and $P \in \textup{Spec}D$:
\begin{eqnarray*}
\tilde I_x &=& \prod_P P^{\tilde r_P(x)}\\
\tilde t_P(x,y)&=& \tilde r_P(x) + \tilde r_{Px}(y) - \tilde r_P(xy) + k_P(x,y),
\end{eqnarray*}
we find
\begin{eqnarray*}
\tilde I_x &=& \sigma_g\left(I_{g^{-1}xg}\right)\alpha(g,g^{-1}xg)\alpha^{-1}(x,g)\\
\tilde r_P(x)&=&r_{Pg}(g^{-1}xg)+k_P(g,g^{-1}xg)-k_P(x,g)\\
\tilde t_P(x,y) &=& t_{Pg}(g^{-1}xg, g^{-1}yg).
\end{eqnarray*}}
\begin{flushleft}\textbf{Proof}\end{flushleft} To make reading easier, we use brackets to do some grouping.
\begin{align*}
u_g T u_g^{-1} &= \bigoplus_{x \in G}u_g I_x u_x u_g^{-1}\\
&= \bigoplus_{x \in G}\sigma_g(I_x)u_g u_x u_g^{-1}\\
&\mathop  = \limits^1  \bigoplus_{x \in G}\sigma_g(I_x)\alpha(g,x)\alpha^{-1}(gxg^{-1},g)u_{g^{-1}xg^{-1}}\\
&\Rightarrow I_x = \sigma_g\left(I_{g^{-1}xg}\right)\alpha(g,g^{-1}xg)\alpha^{-1}(x,g)\\
&\Rightarrow \tilde r_P(x)=r_{Pg}(g^{-1}xg)+k_P(g,g^{-1}xg)-k_P(x,g).
\end{align*}
\begin{align*}\tilde t_P(x,y) =& \tilde r_P(x)+\tilde r_{Px}(y) - \tilde r_P(xy)+k_P(x,y)\\
=&\left[r_{Pg}(g^{-1}xg)+k_P(g,g^{-1}xg)\right]-k_P(x,g)\\
&+\left[r_{Pxg}(g^{-1}yg)+k_{Px}(g,g^{-1}yg)\right]-k_{Px}(y,g)\\
&+\left[-r_{Pg}(g^{-1}xyg)-k_P(g,g^{-1}xyg)\right]+k_P(xy,g)\\
&+k_P(x,y)\\
\mathop  = \limits^2& t_P(g,g^{-1}xg)-r_P(g)+r_P(xg)-k_P(x,g)\\
&+t_{Px}(g,g^{-1}yg)-r_{Px}(g)+r_{Px}(yg)\\
&-t_P(g,g^{-1}xyg)+r_P(g)-r_P(xyg)\\
&+\left[-k_{Px}(y,g)+k_P(xy,g)+k_P(x,y)\right]
\end{align*}
\begin{align*}
\mathop  = \limits^3& t_P(g,g^{-1}xg)+t_{Px}(g,g^{-1}yg)-t_P(g,g^{-1}xyg)\\
&+\left[r_P(xg)-r_{Px}(g)-k_P(x,g)\right]\\
&+\left[r_{Px}(yg)-r_P(xyg)+k_P(x,yg)\right]\\
\mathop  = \limits^{4,5}& t_P(g,g^{-1}xg)-t_P(g,g^{-1}xyg)+r_P(x)-r_P(x)\\
&+\left[t_{Px}(g,g^{-1}yg)-t_P(x,g)+t_P(x,yg)\right]\\
\mathop  = \limits^6& t_P(xg,g^{-1}yg)+t_P(g,g^{-1}xg)-t_P(g,g^{-1}xyg)\\
\mathop  = \limits^7& t_{Pg}(g^{-1}xg,g^{-1}yg). \blok
\end{align*}

Let $A = D \CGR G$ be a crystalline graded ring with $D$ a Dedekind domain.
It is obvious that a maximal graded order will stay maximal after conjugation.  But it happens that the original ring $A$ will not be contained in the conjugated order, see example \ref{ex4}.\\

\subsection{A Special Case}\label{grmax23}

For this section, we take $A = D \CGR G$ be a crystalline graded ring with $D$ a Dedekind domain, but with a special restriction.  Suppose that all prime ideals appearing in the decompostion of the $\alpha(g,h)$, $g,h \in G$ are invariant for the group action, i.e. $Pg=P, \forall g \in G, \forall P$ appearing in the decomposition of the twisted $2$-cocycle.  This is the case if $D$ is a Discrete Valuation Ring, or if all appearing primes $P$ are of the form $(a,b)$ where $a,b \in D^G$.  In this case, all found additive $2$-cocycles and corresponding prime powers as in the sections above are not spectrally twisted.  $\alpha$ however, still is a twisted $2$-cocycle.\\

{\lem With assumptions as above, suppose we have two additive $2$-cocycles $m,m'$ from $G\times G$ to $D \backslash\{0\}$ and suppose $\exists l_P:G \rightarrow \Z$ a map with
\[m'_P(g,h) + l_P(gh)=m_P(g,h)+l_P(g)+l_{Pg}(h),\]
$\forall g,h \in G$.  Furthermore, suppose $m'_P(e,e) = m_P(e,e)=0$, $m_P(g,h) \in \{0,1\}$, $m'_P(g,h) \in \N$, $\forall P \in \textup{Spec}D$, $\forall g,h \in G$, then $l_P(g)\in \N \forall P, \forall g,h \in G.$}\\
\textbf{Proof} Suppose $\exists P, \exists g \in G$ with $l_P(g)<0$ then
\[l_P(g^2)+m'_P(g,g) = m_P(g,g) +l_P(g)+l_P(g) < 0,\]
since $Pg = P$ and $m_P(g,g) \in \{0,1\}$.  We also have
\[l_P(g^3)+m'_P(g,g^2)=m_P(g,g^2)+l_P(g)+l_P(g^2).\]
Sdding these expressions:
\[l_P(g^3) + m'_P(g,g^2)+m'_P(g,g) = m_P(g,g^2)+l_P(g)+l_P(g^2)+m'_P(g,g) <0.\]
Continuing in this fashion we see ($|G|=n$):
\[l_P(g^n)+\sum_{i=1}^{n-1} m'_P(g,g^i) < 0,\]
and since $l_P(g^n) = l_P(e)=0$ this is a contradiction with the fact that $m'_P$ has values in $\N$. $\blok$\\

We can construct from the given ring $A$ the $2$-cocycles $k$ and $m$ as before, but they will not be spectrally twisted.  We construct a maternal order $M$ given a $\gamma$.

{\st With conventions as above, $M$ is the unique maximal graded order over $D$ containing $A$.}\\
\textbf{Proof} Let $T$ be a graded $D$-order in $K \CGR G$, where $K$ is the quotient field of $D$.  We can associate to $T$ the maps $t$ and $r$ like in Section \ref{grmax8}, i.e.
\[T = \bigoplus J_g u_g, \ \ \ J_g = \prod_{P \in \textup{Spec}D} P^{r_P(g)},\]
\[t_P(g,h)= -r_P(gh)+r_P(g)+r_P(h)+k_P(x,y) \geq 0.\]
This implies, since
\[m_P(g,h)=-a_P(gh)+a_P(g)+a_P(h)+k_P(g,h),\]
that
\[t_P(g,h)= -r_P(gh)+a_P(gh)+r_P(g)-a_P(g)+r_P(h)-a_P(h)+m_P(x,y) \geq 0,\]
$\forall P \in \textup{Spec}D$ and $\forall g,h \in G$.  By the previous lemma we find $r_P(g)-a_P(g)\geq 0$, $\forall P \in \textup{Spec}D$, $\forall g \in G$.  In other words, $\forall g \in G$:
\[J_g = \prod_{P \in \textup{Spec}D} P^{r_P(g)}\subset \prod_{P \in \textup{Spec}D} P^{a_P(g)}=I_g \Rightarrow T \subset A.\]
The theorem now follows. $\blok$\\

\section{Examples}

\section{Example 1: $1$ maximal graded order}\label{ex2}
Consider $\Z[i] \CGR \Z_4$ where
\[\sigma: \Z_4 \rightarrow \textup{Aut}\Z[i]:\bar 0, \bar 2 \mapsto \textup{Id}, \ \ , \bar 1, \bar 3 \mapsto \bar \cdot,\]
where $\bar \cdot$ is the standard conjugation in $\C$.  We set $\alpha$:
\[
\begin{array}{*{20}c}
   \alpha  &\vline &  \bar 0 & \bar1 & \bar2 & \bar3  \\
\hline
   \bar0 &\vline &  1 & 1 & 1 & 1  \\
   \bar1 &\vline &  1 & 2 & {1 - 2i} & {2 - 4i}  \\
   \bar2 &\vline &  1 & {1 - 2i} & 5 & {1 + 2i}  \\
   \bar3 &\vline &  1 & {2 + 4i} & {1 + 2i} & 2  \\
 \end{array}
\]
The relevant primes are
\begin{eqnarray*}
P_1 &=& 1+2i,\\
P_2 &=& 1-2i,\\
P_3 &=& 1+i.
\end{eqnarray*}
We calculate $k_{P_1},k_{P_2},k_{P_3}$:
\[
\begin{array}{*{20}c}
   k_{P_1}  &\vline &  \bar 0 & \bar1 & \bar2 & \bar3  \\
\hline
   \bar0 &\vline &  0 & 0 & 0 & 0  \\
   \bar1 &\vline &  0 & 0 & 0 & 0  \\
   \bar2 &\vline &  0 & 0 & 1 & 1  \\
   \bar3 &\vline &  0 & 1 & 1 & 0  \\
 \end{array}
\quad \quad \quad\begin{array}{*{20}c}
   k_{P_2}  &\vline &  \bar 0 & \bar1 & \bar2 & \bar3  \\
\hline
   \bar0 &\vline &  0 & 0 & 0 & 0  \\
   \bar1 &\vline &  0 & 0 & 1 & 1  \\
   \bar2 &\vline &  0 & 1 & 1 & 0  \\
   \bar3 &\vline &  0 & 0 & 0 & 0  \\
 \end{array}
\quad \quad \quad
\begin{array}{*{20}c}
   k_{P_3}  &\vline &  \bar 0 & \bar1 & \bar2 & \bar3  \\
\hline
   \bar0 &\vline &  0 & 0 & 0 & 0  \\
   \bar1 &\vline &  0 & 2 & 0 & 2  \\
   \bar2 &\vline &  0 & 0 & 0 & 0  \\
   \bar3 &\vline &  0 & 2 & 0 & 2  \\
 \end{array}
\]
We calculate $\mu$ such that $\alpha^4(x,y)=\mu(xy)^{-1}\mu(x)\sigma_x(\mu(y))$.  We have
\begin{eqnarray*}
\alpha^{4}(\bar2,\bar2) &=& \mu(\bar2)^2 = 625,\\
\alpha^{4}(\bar1,\bar3) &=& \mu(\bar1)\overline{\mu(\bar3)}=2^4(1-2i)^4,\\
\alpha^{4}(\bar1,\bar3) &=& \overline{\mu(\bar1)}\mu(\bar3)=2^4(1+2i)^4,\\
\alpha^4(\bar1,\bar1)&=&\mu(\bar2)^{-1}\mu(\bar1)\overline{\mu(\bar1)}= \frac{1}{25}|\mu(\bar1)|^2 = 16 \Rightarrow |\mu(\bar1)|^2 = 400.
\end{eqnarray*}
From which we can conclude
\begin{eqnarray*}
\mu(\bar1)&=& 4 (1-2i)^2,\\
\mu(\bar2)&=& 25,\\
\mu(\bar3)&=& 4 (1+2i)^2.
\end{eqnarray*}
We now calculate $\gamma$ and $1/4 \gamma$:
\[
\begin{array}{*{20}c}
   \gamma  &\vline &  {\bar 0} & {\bar 1} & {\bar 2} & {\bar 3}  \\
\hline
   {P_1 } &\vline &  0 & 0 & 2 & 2  \\
   {P_2 } &\vline &  0 & 2 & 2 & 0  \\
   {P_3 } &\vline &  0 & 4 & 0 & 4  \\

 \end{array} \quad \quad \quad \begin{array}{*{20}c}
   {\frac{1}
{4}\gamma } &\vline &  {\bar 0} & {\bar 1} & {\bar 2} & {\bar 3}  \\
\hline
   {P_1 } &\vline &  0 & 0 & {\frac{1}
{2}} & {\frac{1}
{2}}  \\
   {P_2 } &\vline &  0 & {\frac{1}
{2}} & {\frac{1}
{2}} & 0  \\
   {P_3 } &\vline &  0 & 1 & 0 & 1  \\

 \end{array}
\]
And then $a$ ($0 \leq a_{P_j} + 1/4 \gamma_{P_j}<1$):
\[
\begin{array}{*{20}c}
   a &\vline &  {\bar 0} & {\bar 1} & {\bar 2} & {\bar 3}  \\
\hline
  {P_1 } &\vline &  0 & 0 & 0 & 0  \\
   {P_2 } &\vline &  0 & 0 & 0 & 0  \\
   {P_3 } &\vline &  0 & { - 1} & 0 & { - 1}  \\

 \end{array}
\]
So the maternal order is
\[\Z[i]u_{\bar 0} + (1+i)^{-1}\Z[i]u_{\bar 1}+\Z[i]u_{\bar 2} + (1+2)^{-1}\Z[i]u_{\bar 3},\]
We check for maximality:
\begin{eqnarray*}
m_{P_1}(\bar 3, \bar 1)&=&1,\\
m_{P_1}(\bar 2, \bar 2)&=&1,\\
m_{P_2}(\bar 1, \bar 3)&=&1,\\
m_{P_2}(\bar 2, \bar 2)&=&1.
\end{eqnarray*}
\[t_{P_1}(\bar2,\bar2) = r_{P_1}(\bar2)+r_{P_1}(\bar2)+k_{P_1}(\bar 2, \bar 2),\]
and so $r_{P_1}(\bar 2) = a_{P_1}(\bar 2)$ and similarly $r_{P_2}(\bar 2) = a_{P_2}(\bar 2)$.  We now check case (\ref{grmax11}) for $x = \bar 3$, $P=P_1$.
\[
\begin{array}{*{20}c}
   {m_{P_1 } (\bar 3,\bar 1)} \hfill & {x = \bar 3} \hfill & {y = \bar 1} \hfill & {m_{P_2 } (\bar 1,\bar 3) = 1} \hfill  \\
   {} \hfill & {} \hfill & {y = \bar 2} \hfill & {m_{P_1 } (\bar 2,\bar 3) = 1} \hfill  \\
   {} \hfill & {} \hfill & {y = \bar 3} \hfill & {m_{P_2 } (\bar 3,\bar 3) = 0} \hfill  \\
 \end{array}
\]
This means we need to modify $r_{P_2}(\bar 3 + \bar 3)=r_{P_2}(\bar 2)$ which is impossible.  The same problem appears in the case of $m_{P_2}(\bar 1,\bar 3)$.  This means that the maternal order itself is maximal.

\section{Example 2: $2$ maximal graded orders}\label{ex3}
Consider $\Z[i] \CGR \Z_2$ where
\[\sigma: \Z_2 \rightarrow \textup{Aut}\Z[i]: \bar 0 \mapsto \textup{Id},\ \ ,\bar1 \mapsto \bar \cdot.\]
Put $\alpha(\bar 1, \bar 1)=5$, then the relevant primes are $P_1 = 1+2i$, $P_2 = 1-2i$.  Calculating the $k$ yields $k_{P_1}(\bar 1, \bar 1) = k_{P_2}(\bar 1, \bar 1) =1$.  $\alpha^2(\bar 1, \bar 1) = 25$ and so $\mu(\bar 1)=5$.  We now find that
\[
\gamma_{P_1}(\bar1) = 1, \ \gamma_{P_2}(\bar1)=1, \ a_{P_1}(\bar1) = 0, \ a_{P_2}(\bar1)=0, \ m_{P_1}(\bar1,\bar1) = 1, \ m_{P_2}(\bar1, \bar1)=1, \ \ 
\]
This means that the maternal order is
\[\Z[i]u_{\bar0}+\Z[i]u_{\bar1}.\]
We can either
\begin{enumerate}
	\item \begin{eqnarray*}
	r_{P_1}(\bar 1)&=& a_{P_1}(\bar 1)-1 = -1\\
	r_{P_2}(\bar 1)&=& a_{P_2}(\bar 1) = 0,
	\end{eqnarray*}
	\item \begin{eqnarray*}
	r_{P_2}(\bar 1)&=& a_{P_2}(\bar 1)-1 = -1\\
	r_{P_1}(\bar 1)&=& a_{P_1}(\bar 1) = 0.
	\end{eqnarray*}
	
\end{enumerate}
And we get the following maximal graded orders
\begin{eqnarray*}
T_1 &=&\Z[i]u_{\bar 0}+(1+2i)^{-1}\Z[i]u_{\bar 1},\\
T_2 &=&\Z[i]u_{\bar 0}+(1-2i)^{-1}\Z[i]u_{\bar 1}.
\end{eqnarray*}
As one can check, both are even strongly graded.

\section{Example 3: $2$ maximal graded orders}\label{ex4}

Consider $\Z[i] \CGR \Z_4$ where
\[\sigma: \Z_4 \rightarrow \textup{Aut}\Z[i]:\bar 0, \bar 2 \mapsto \textup{Id}, \ \ , \bar 1, \bar 3 \mapsto \bar \cdot,\]
where $\bar \cdot$ is the standard complex conjugation.  We set $\alpha$:
\[
\begin{array}{*{20}c}
   \alpha  &\vline &  0 & 1 & 2 & 3  \\
\hline
   0 &\vline &  1 & 1 & 1 & 1  \\
   1 &\vline &  1 & 4 & {(1 + i)(1 + 2i)} & {20(1 + i)(1 + 2i)}  \\
   2 &\vline &  1 & {(1 + i)(1 + 2i)} & {50} & {5(1 - i)(1 - 2i)}  \\
   3 &\vline &  1 & {20(1 - i)(1 - 2i)} & {5(1 - i)(1 - 2i)} & {20}  \\

 \end{array}
\]
The relevant primes are
\begin{eqnarray*}
P_1 &=& 1+2i,\\
P_2 &=& 1-2i,\\
P_3 &=& 1+i.
\end{eqnarray*}
We calculate $k_{P_1},k_{P_2},k_{P_3}$:
\[
\begin{array}{*{20}c}
   {k_{P_1 } } &\vline &  0 & 1 & 2 & 3  \\
\hline
   0 &\vline &  0 & 0 & 0 & 0  \\
   1 &\vline &  0 & 0 & 1 & 2  \\
   2 &\vline &  0 & 1 & 2 & 1  \\
   3 &\vline &  0 & 1 & 1 & 1  \\

 \end{array}
\quad \quad \quad \begin{array}{*{20}c}
   {k_{P_2 } } &\vline &  0 & 1 & 2 & 3  \\
\hline
   0 &\vline &  0 & 0 & 0 & 0  \\
   1 &\vline &  0 & 0 & 0 & 1  \\
   2 &\vline &  0 & 0 & 2 & 2  \\
   3 &\vline &  0 & 2 & 2 & 1  \\

 \end{array}
\quad \quad \quad\begin{array}{*{20}c}
   {k_{P_3 } } &\vline &  0 & 1 & 2 & 3  \\
\hline
   0 &\vline &  0 & 0 & 0 & 0  \\
   1 &\vline &  0 & 4 & 1 & 5  \\
   2 &\vline &  0 & 1 & 2 & 1  \\
   3 &\vline &  0 & 5 & 1 & 4  \\

 \end{array}
\]
From this we can calculate the possibilities for $\gamma$:
\[
\begin{array}{*{20}c}
   \gamma  &\vline &  0 & 1 & 2 & 3  \\
\hline
   {P_1 } &\vline &  0 & s & 4 & s  \\
   {P_2 } &\vline &  0 & {4 - s} & 4 & {8 - s}  \\
   {P_3 } &\vline &  0 & {10} & 4 & {10}  \\

 \end{array}
\]
where $s \in \Z$ and $0 \leq s \leq 4$.  And so we can find 3 different maternal orders, namely one corresponding to $s = 4$, call it $M^A$, one corresponding to $s = 0$, $M^B$ and one where $1 \leq s \leq 3$, $L$.\\
\begin{eqnarray*}
M^A &=& \Z[i]u_0 + (1+2i)^{-1}(1+i)^{-2}u_1\\
&+& (1+2i)^{-1}(1-2i)^{-1}(1+i)^{-1}u_2 + (1+2i)^{-1}(1-2i)^{-1}(1+i)^{-2}u_3.\\
&&\\
M^B&=&\Z[i]u_0 + (1-2i)^{-1}(1+i)^{-2}u_1\\
&+&(1-2i)^{-1}(1+2i)^{-1}(1+i)^{-1}u_2 + (1-2i)^{-2}(1+i)^{-2}u_3.\\
&&\\
L &=& \Z[i]u_0 + (1+i)^{-2}u_1\\
&+& (1-2i)^{-1}(1+i)^{-1}(1+2i)^{-1}u_2 + (1-2i)^{-1}(1+i)^{-2}u_3.
\end{eqnarray*}

Calculation and the method described in section \ref{grmax8} reveals that $M_1$ and $M_2$ are maximal graded orders, and $L$ is not.  Depending on the choices made in the process of finding orders above $L$, we find $M_1$ and $M_2$.\\

Looking at conjugation, we find the following:
\begin{eqnarray*}
\Psi_3(M^A)&=&M^A,\\
\Psi_1(M^A)&=&M^B,\\
\Psi_1 \circ \Psi_1 &=& \textup{Id},\\
\Psi_3 \circ \Psi_3 &=& \textup{Id}.
\end{eqnarray*}
We also find that $\Psi_3\left(\Psi_1(M^A)\right)$ does not contain $A$ anymore!  Putting it in a picture:

\begin{figure}[!ht]
  \begin{center}
    \includegraphics[width=7.5cm]{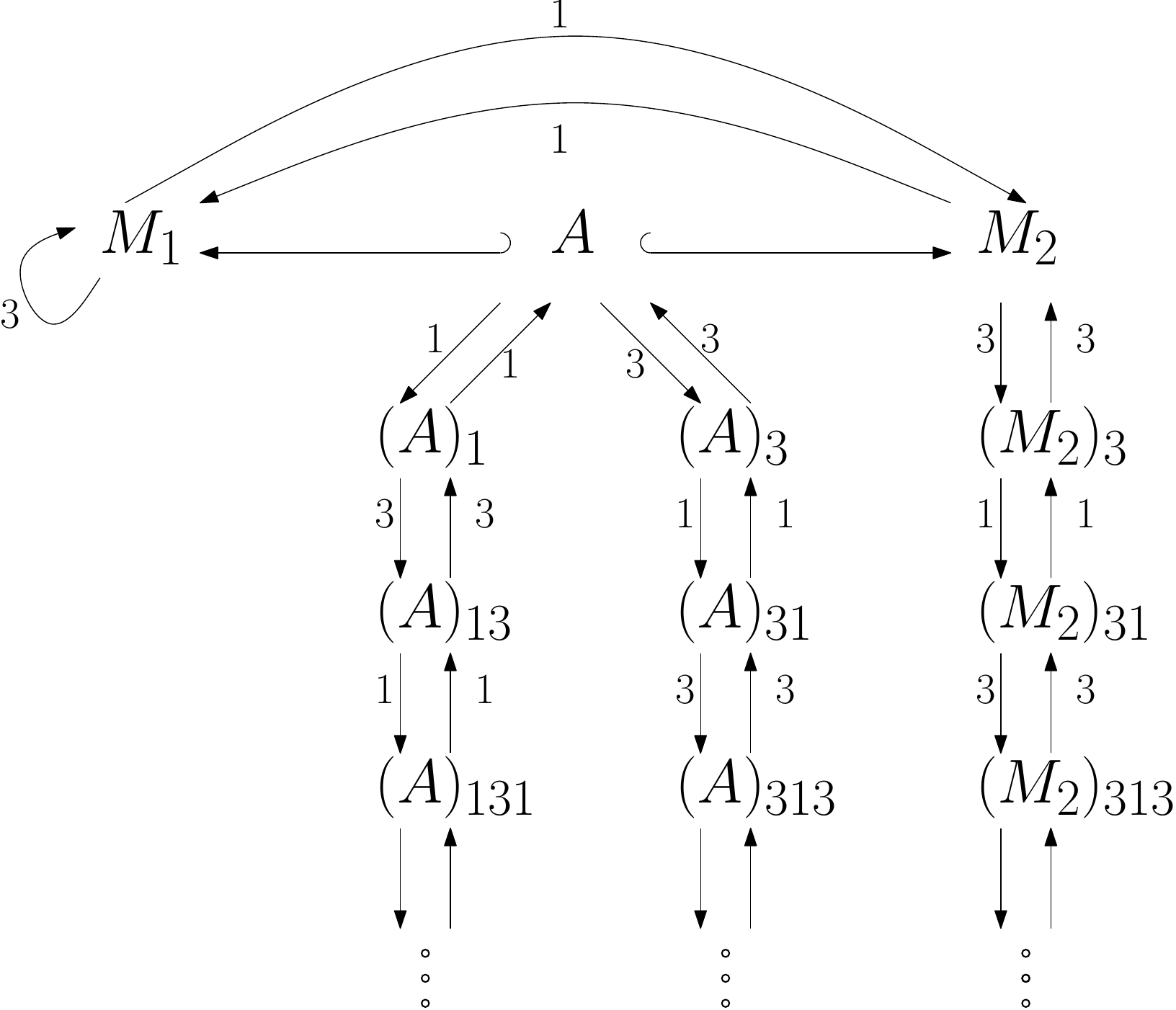}
  \end{center}
\end{figure}

\end{document}